\documentclass{amsart}
\usepackage{amsfonts,latexsym,amssymb,amsmath,amsthm}
\usepackage{enumerate}
\usepackage{captdef}
\usepackage[dvips]{graphicx} 
\usepackage{enumerate,graphicx,graphpap}
\usepackage[matrix]{xy}
\usepackage{cases}
\usepackage{epsfig}
\usepackage{color}
\input xy
\xyoption{all} \xyoption{poly}

\theoremstyle{plain}
\newtheorem{theorem}{Theorem}[section]
\newtheorem{lema}{Lemma}[section]

\newtheorem{remark}[theorem]{Remark}

\theoremstyle{definition}

\newtheoremstyle{teoremacita}
{3pt}
{3pt}
{\itshape}
{}
{\bfseries}
{}
{ }
{\thmname{#1}\thmnumber{ #2'}\thmnote{ #3}.}
\theoremstyle{teoremacita} \newtheorem*{teor*}{}

\newcommand{\be}{\begin{enumerate}}
\newcommand{\ee}{\end{enumerate}}

\newcommand{\bi}{\begin{itemize}}
\newcommand{\ei}{\end{itemize}}

\def\N{\mathbb N}

\def\R{\mathbb R}
\def\C{\mathbb C}

\def\a{\alpha}

\def\d{\partial}

\begin{document}

\title{Can we detect Gaussian curvature by counting paths and measuring their length?}

\author{Leonardo A. Cano G.}
\address[Leonardo Cano]{Departamento de Matem\'{a}ticas\\
Facultad de Ciencias, Universidad Nacional de Colombia\\
Carrera 45 26-85\\
Bogot\'{a} - Colombia}
\email{lcanog@unal.edu.co}

\author{Sergio A. Carrillo}
\address[Sergio A. Carrillo]{Programa de Matem\'{a}ticas\\
Escuela de Ciencias Exactas e Ingenier\'{i}a, Universidad Sergio Arboleda\\
Calle 74 14-14\\
Bogot\'{a} - Colombia}
\email{sergio.carrillo@usa.edu.co}


\thanks{The second author was partially supported by the Ministerio de Econom\'{\i}a y Competitividad from Spain, under the Project ``\'{A}lgebra y geometr\'{\i}a en sistemas din\'{a}micos y foliaciones singulares" (Ref.: MTM2016-77642-C2-1-P)}

\begin{abstract} The aim of this paper is to associate a measure for certain sets of paths in the Euclidean plane $\R^2$ with fixed starting and ending points. Then, working on parameterized surfaces with a specific Riemannian metric, we define and calculate the integral of the length over the set of paths obtained as the image of the initial paths in $\R^2$ under the given parameterization. Moreover, we prove that this integral is given by the average of the lengths of the external paths times the measure of the set of paths if and only if the surface has Gaussian curvature equal to zero.
\end{abstract}

\maketitle

\

\section{Introduction}

In the theory of smooth surfaces the Gaussian curvature plays a fundamental role as a tool to measure how these objects bend in Euclidean space. It is an intrinsic function depending only on distances that are measured on the surface, but not in how the surface is isometrically embedded in $\R^3$. This result is the famous Gauss' Theorema Egregium.

Can we detect this curvature by counting paths and measuring their length in a given surface? To fix ideas, let us consider a connected surface $S$ and two given points $p,q\in S$. We will not consider all possible paths on $S$ from $p$ to $q$, but only simple ones as we shall explain. To settle down the problem, we start with the simplest case, i.e., the Euclidean plane $\R^2$, of constant curvature equal to zero, when equipped with Cartesian coordinates $(x,y)$ and the Riemannian metric $dx^2+dy^2$. The idea relies on looking only at paths built with straight lines parallel to the axes, and only moving forward. In other words, we concatenate a finite number of pieces of the flows $((x,y),s)\mapsto (x+s,y)$, $((x,y),s)\mapsto (x,y+s)$, corresponding to the orthogonal vector fields $\d_x$ and $\d_y$ coming from the given coordinates. Then we can reach $q=(x_1,y_1)$ from $p=(x_0,y_0)$ if and only if $x_1\geq x_0$, $y_1\geq y_0$ and the time we need to do this is given by $t=x_1-x_0+y_1-y_0$. In this way we have established the set $\Gamma_{p,q}(t)$ of all admissible paths. The first question is: how can we count or give a measure to this set? In this particular case, every element of $\Gamma_{p,q}(t)$ has the same length, i.e., $t$. Thus, independently of how we define the measure value  $m(\Gamma_{p,q}(t))$, when we add all these lengths, the result will be $t\cdot m(\Gamma_{p,q}(t))$. The conclusion is the following: \textit{in curvature zero the ``sum'' of all these lengths is equal to the average of the lengths of the external paths times the ``number'' of paths}. The external paths are the one from $(x_0,y_0)$ to $(x_1,y_0)$ and then to $(x_1,y_1)$, and the other from $(x_0,y_0)$ to $(x_0,y_1)$ and then to $(x_1,y_1)$.

The aim of this note is to show that we can effectively define  $m(\Gamma_{p,q}(t))$ and the ``sum'', or better, the integral of the length over $\Gamma_{p,q}(t)$, to prove that the phenomenon explained in the previous paragraph only occurs in the case of surfaces of curvature zero. We will transfer the situation of the Euclidean plane to $S$ by using a local  parameterization $\Psi:U\subseteq\R^2\rightarrow S$ containing $p=\Psi(x_0,y_0)$ and $q=\Psi(x_1,y_1)$. Our calculations work for metrics of the form  \begin{equation}\label{Metric General} h'(x)^2dx^2+f'(x)^2dy^2,
\end{equation} where $h, f:\R\rightarrow \R$ are strictly increasing smooth functions. We remark that, when $|f''(s)|< h'(s)$, surfaces of revolution realize these spaces (see Section \ref{Sec: Adding all the lengths} for details). Since the classical surfaces of constant curvature, i.e, the plane $\R^2$, the sphere $S^2$ and the Poincar\'{e} half plane $\mathbb{H}=\R_{>0}\times \R$, can be realized locally as open revolution surfaces in $\R^3$, our result is valid for these surfaces as well.  

We will see that a good definition for the measure of the set of paths $\Gamma_{p,q}(t)$ explained above is given by   $m(\Gamma_{p,q}(t))={t\brace x_1-x_0}={t\brace y_1-y_0}$, where \begin{equation}\label{Binomial Coefficients}
{t\brace a}:=2\sum_{n=0}^\infty \frac{a^n(t-a)^n}{n!^2}+t\sum_{n=0}^\infty\frac{a^{n}(t-a)^n}{(n+1)!n!},
\end{equation} denotes the so called  \textit{continuous binomial coefficient} defined in \cite{CD}. The name and suggestive notation come from the analogy with the binomial coefficient $\binom{n+m}{n}$ which counts the paths in the lattice $\N^2\subset\R^2$, from $(0,0)$ to $(m,n)\in\N^2$. Then, our main result is an explicit formula for the integral of the length associated to the Riemannian metric (\ref{Metric General}), which is the content of:

\begin{theorem}\label{Th. Main} The integral of the length over all admissible paths $\Gamma_{p,q}(t)$ from $p=\Psi(x_0, y_0)$ to $q=\Psi(x_1, y_1)$ in time $t>0$ is given by  \begin{align}\label{Main formula} \int\limits_{\Gamma_{p,q}(t)} l(\gamma)d(\gamma)=&\left(h(x_1)-h(x_0)+\frac{(y_1-y_0)}{2}\left(f'(x_1)+f'(x_0)\right)\right)F(x_1-x_0,y_1-y_0)+\\ \nonumber
	&\left(f(x_1)-f(x_0)-\frac{(x_1-x_0)}{2}\left(f'(x_1)+f'(x_0)\right)\right)\frac{\d^2 F}{\d \tau^2}(x_1-x_0,y_1-y_0),
	\end{align} where $t=x_1-x_0+y_1-y_0$ and  $F(\tau,s):={\tau+s\brace s}$. In particular, the integral is given by the measure of $\Gamma_{p,q}(t)$ times the average of the length of the external paths if and only if $f$ is a quadratic function. For the surfaces of constant curvature, this happens only if the curvature is zero.
\end{theorem}

\

This formula gives us a new insight on how the curvature of a surface affects the length of paths defined over it.  We will prove it by studying a quite simple definition of a path integral and looking at the geometric information it contains.

It is worth remarking that several approaches to path integrals exist in the literature. The most famous ones are Feynman and Wiener integrals. Feynman integrals, although not fully formalized from a mathematical point of view, are used in Quantum Mechanics to capture physical information, see \cite{Feynman,CD2001}. The Wiener integrals are defined on sets of continuous paths in Euclidean space and model Brownian motion. However, the differential paths have Wiener measure zero, see e.g., \cite[Thm. 1.1]{CD2001}, and hence do not have a notion of length. It could be interesting to see if the integral of length makes some sense in any of the mathematical formulations of Feynman integrals.

\


\section{Paths as simplexes and their measure}\label{Sec: Paths as simplexes and their measure}

Coming back to the case of the plane $\R^2$, a path in $\Gamma_{p,q}(t)$ is determined, first, by the order in which we glue a finite number of segments given by the flows of $\d_x$ and $\d_y$, and second, by the length of these segments: for the horizontal (resp. vertical) part, the length is $x_1-x_0$, (resp. $y_1-y_0$). Such order can be represented by a finite sequence $\textbf{c}$ of the numbers $1$ and $2$, where $1$ represents an horizontal segment and $2$ represents a vertical one, i.e., by an element of the set $$C(n):= \{\textbf{c}=(c_0, c_1,\dots,c_n)\in \{1,2\}^{n+1} \,|\, c_{i-1}\neq c_i, i=1,\dots,n. \},\quad \text{ for some }n\in\N^\ast.$$ In this way, we can write \begin{equation}\label{Gamma pq}
\Gamma_{p,q}(t)=\bigsqcup_{n\in\N^\ast} \bigsqcup_{\textbf{c}\in C(n)} \Gamma_{p,q}^{\textbf{c}}(t),
\end{equation} as a disjoint union, where each $\Gamma_{p,q}^{\textbf{c}}(t)$ denotes the set of paths built with order $\textbf{c}$. If $\textbf{c}\in C(n)$, a path $\gamma\in \Gamma_{p,q}^{\textbf{c}}(t)$ corresponds to a vector $\textbf{s}=(s_0,s_1,\dots,s_n)\in\R^{n+1}_{>0}$, where $s_j$ is the length of the $j$th segment of $\gamma$. We can distinguish four cases:

\begin{enumerate}
	\item[I.1.]  If $n=2m-1$ and $\textbf{c}=(1,2,\dots,1,2)$, $\textbf{s}$ satisfies $\sum_{j=0}^{m-1} s_{2j}=x_1-x_0$ and $\sum_{j=0}^{m-1} s_{2j+1}=y_1-y_0$. 
	
	\item[I.2.]  If $n=2m-1$ and $\textbf{c}=(2,1,\dots,2,1)$, $\textbf{s}$ satisfies  $\sum_{j=0}^{m-1} s_{2j} =y_1-y_0$ and $\sum_{j=0}^{m-1} s_{2j+1}=x_1-x_0$. 
	
	\item[II.1.] If $n=2m$ and  $\textbf{c}=(1,2,\dots,1,2,1)$, $\textbf{s}$ satisfies  $\sum_{j=0}^{m} s_{2j}=x_1-x_0$ and $\sum_{j=0}^{m-1} s_{2j+1}=y_1-y_0$. 
	\item[II.2.] If $n=2m$ and  $\textbf{c}=(2,1,\dots,2,1,2)$, $\textbf{s}$ satisfies $\sum_{j=0}^{m} s_{2m}=y_1-y_0$ and $\sum_{j=0}^{m-1} s_{2j+1}=x_1-x_0$.
\end{enumerate} 

In all cases, $\Gamma_{p,q}^{\textbf{c}}(t)$ is given, up to a permutation of the variables, by the Cartesian product of two simplexes. Thus we can reduce the problem to measure usual simplexes in Euclidean space.

Given $\tau\geq 0$, we denote by $$\Delta_n^\tau:=\{\textbf{s}=(s_0,s_1,\dots,s_n)\in\R_{\geq0}^{n+1} \,|\, s_0+s_1+\cdots+s_n=\tau\},$$ the $n$--symplex generated by the canonical vectors of $\R^{n+1}$, multiplied by $\tau$. We can parameterize it by using
\begin{equation}\label{chart varphi}\varphi_n^\tau:\R^n\longrightarrow \R^{n+1},\quad 
\varphi(l_1,l_2,\dots,l_n)=(l_1,l_2-l_1,\dots,l_n-l_{n-1},\tau-l_n),
\end{equation} restricted to $W_n^\tau:=\{(l_1,\dots,l_n)\in\R^n \,|\,  0\leq l_1\leq\dots\leq l_n\leq \tau\}$. The way we measure $\Delta_n^\tau$ is by calculating the $n$--volume of $W_n^\tau$, i.e., by declaring that  $$m(\Delta_n^\tau):=\int\limits_0^\tau\int\limits_0^{l_{n}}\cdots \int\limits_0^{l_2} 1\, dl_1\cdots dl_n=\frac{\tau^n}{n!}.$$ If $n=0$, we interpret this value as $1$. Note that the variables involved are related by the formulas  \begin{equation}\label{change of varialbles l}
l_1=s_0,\quad l_2=s_0+s_1,\quad \dots,\quad l_{n}=s_0+s_1+\cdots+s_{n-1}.\end{equation} 

By using the parameterization $\varphi_{n_1}^{\tau_1}\times \varphi_{n_2}^{\tau_2}:W_{n_1}^{\tau_1}\times W_{n_2}^{\tau_2}\rightarrow \Delta_{n_1}^{\tau_1}\times \Delta_{n_2}^{\tau_2}$, we find that $$m(\Delta_{n_1}^{\tau_1}\times \Delta_{n_2}^{\tau_2})=m(\Delta_{n_1}^{\tau_1})\cdot m(\Delta_{n_2}^{\tau_2}),$$ formula that can be interpreted as a multiplicative counting principle in our setting. Moreover,  if we fix $0\leq m<n$, $\Delta_n^\tau$ can be written as the disjoint union $\bigsqcup_{0\leq s\leq \tau}\Delta_m^s\times \Delta_{n-1-m}^{\tau-s}$. The measure we have introduced is compatible with this decomposition in the following sense: $$\int_0^\tau m\left(\Delta_m^s\times \Delta_{n-1-m}^{\tau-s}\right)ds=\int_0^\tau \frac{s^m}{m!}\frac{(\tau-s)^{n-1-m}}{(n-1-m)!}ds=\frac{\tau^n}{n!}=m(\Delta_n^\tau).$$ Note we have used here the elementary formula for the Beta function  $\int_0^1 u^p(1-u)^qdu=\frac{p!q!}{(p+q+1)!}$, $p,q\in\N$.

In this order of ideas, we have that  $m(\Gamma^{\textbf{c}}_{p,q}(t))=\frac{(x_1-x_0)^{m-1}}{(m-1)!}\frac{(y_1-y_0)^{m-1}}{(m-1)!}$ for the cases I.1 and I.2, $m(\Gamma^{\textbf{c}}_{p,q}(t))=\frac{(x_1-x_0)^{m}}{m!}\frac{(y_1-y_0)^{m-1}}{(m-1)!}$ for the case II.1, and $m(\Gamma^{\textbf{c}}_{p,q}(t))=\frac{(x_1-x_0)^{m-1}}{(m-1)!}\frac{(y_1-y_0)^m}{m!}$ for the case II.2. Finally, the measure of $\Gamma_{p,q}(t)$ is defined as the sum of the measures of its parts, \begin{equation}m\left(\Gamma_{p,q}(t)\right):=\sum_{n=1}^\infty \sum_{\textbf{c}\in C(n)} m(\Gamma_{p,q}^\textbf{c}(t)), \end{equation} in agreement with the definition of the binomial coefficient ${t\brace x_1-x_0}$ given in the Introduction.

\begin{remark} Our approach to associate this measure is based on the work \cite{CD1}. There the authors work more generally with a manifold $M$ and a finite number of vector fields $X_1,\dots, X_k$ defined over it. In this framework, if $p,q\in M$ and $t>0$ are given, the set of admissible paths $\Gamma_{p,q}(t)$  consists of piece-wise smooth paths from $p$ to $q$ formed by concatenations of the flows $\phi_j$ of the $X_j$ using a total time $t$. Here there is also a decomposition of $\Gamma_{p,q}(t)$ analogous to (\ref{Gamma pq}) and if $\mathbf{c}=(c_0.c_1,\dots,c_n)$ is a given order, 
	$\Gamma_{p,q}^{\mathbf{c}}(t)$ can be embedded in $\Delta_n^t$.
\end{remark}

\section{On the continuous binomial coefficients}\label{Sec: Why the continuous binomial coefficients?}

The motivation to introduce the power series ${t\brace a}$, which converges for every $t,a\in\C$, came from measure or ``count'' elements of $\Gamma_{p,q}(t)$, thus we regard it as a continuous analogue (not extension) of the classical binomial coefficients. For extensions to complex variables, it is enough to write the usual binomial coefficients in terms of the Gamma function. For a more detailed study of this extension to real variables, see e.g. \cite{Salwinski}.

To show the strong parallel with the discrete binomial coefficients, we highlight the following two properties $${\tau+s\brace s}={\tau+s\brace \tau},\quad \frac{\d^2}{\d \tau\d s}{\tau+s\brace s}={\tau+s\brace s},$$ which corresponds to the discrete versions  $\binom{n+m}{n}=\binom{n+m}{m}$ and $D_nD_m\binom{n+m}{m}=\binom{n+m}{m}$, respectively, and $D_nf=f(n+1)-f(n)$ is a difference operator. There are more analogies, for instance, continuous versions of Chu--Vandermonde's formula, see \cite{W2}.

To perform some calculations we will need later, we can use the \textit{Bessel--Clifford function of the first kind} \cite[p. 358]{AS} which are defined by the power series \begin{equation*}\label{Def. Bessel Cliford Function} C_\nu(z):=\sum_{m=0}^\infty \frac{z^m}{m!\Gamma(1+m+\nu)},\quad  z, \nu\in\C,
\end{equation*} and where $\Gamma$ denotes the classical gamma function. From the very definition, we note that $\frac{dC_\nu}{dz}=C_{\nu+1}$ and also \begin{equation}\label{Cn+2 Cn+1 Cn} zC_{\nu+2}(z)+(\nu+1)C_{\nu+1}(z)=C_\nu(z).\end{equation}

We can use the functions $C_0$ and $C_1$ to write \begin{equation}\label{BC in terms BC}
F(\tau,s):={\tau+s\brace s}=2C_0(s\tau)+(\tau+s)C_1(s\tau),\end{equation} and also the remaining $C_n$ to write the derivatives \begin{equation}\label{Formulas derivatives}
\frac{\d^n F}{\d \tau^n}(\tau,s)=s^{n-1}(2s+n)C_{n}(s\tau)+(\tau+s)s^nC_{n+1}(s\tau)=\frac{\d^n F}{\d s^n}(s,\tau),\quad n\geq1.
\end{equation} These equations are easily proved by using induction on $n$. For later use we need the following two formulas: \begin{align}
\label{Formula 1} \frac{\d^2 F}{\d \tau^2}(\tau,s)&=\sum_{m=0}^\infty 2\frac{\tau^m}{m!}\frac{s^{m+2}}{(m+2)!}+\frac{\tau^{m}}{m!}\frac{s^{m+3}}{(m+3)!}+\frac{\tau^{m}}{m!}\frac{s^{m+1}}{(m+1)!},\\
\label{Formula 2} sF(\tau,s)-\tau\frac{\d^2 F}{\d \tau^2}(\tau,s)&=2\sum_{m=0}^\infty \frac{s^{m+1}}{(m+1)!}\frac{\tau^m}{m!}+\frac{s^{m+2}}{(m+2)!}\frac{\tau^{m}}{m}.
\end{align}

To prove (\ref{Formula 1}), simply use (\ref{Formulas derivatives}) for $n=2$, and (\ref{Cn+2 Cn+1 Cn}) for $\nu=1$ to obtain $\frac{\d^2 F}{\d \tau^2}(\tau,s)=2s(s+1)C_2(s\tau)+(\tau+s)s^2 C_3(s\tau)=2s^2C_2(s\tau)+s^3C_3(s\tau)+s(2C_2(s\tau)+s\tau C_3(s\tau))=2s^2C_2(s\tau)+s^3C_3(s\tau)+sC_1(s\tau)$. In the same way, we use (\ref{Cn+2 Cn+1 Cn}) for $\nu=0$ and $\nu=1$ so see that the left-hand side of (\ref{Formula 2}) is in fact \begin{align*}
&s\left(2C_0(s\tau)+(\tau+s)C_1(s\tau)\right)-\tau\left(2s(s+1)C_2(s\tau)+s^2(\tau+s)C_3(s\tau)\right)\\
=&s\left(2C_0(s\tau)+(\tau+s)C_1(s\tau)-2\tau(s+1)C_2(s\tau)-s\tau(\tau+s)C_3(s\tau)\right)\\
=&s\left(2\left[s\tau C_2(s\tau)+C_1(s\tau)\right]+(\tau+s)C_1(s\tau)-2\tau(s+1)C_2(s\tau)-(\tau+s)\left[C_1(s\tau)-2C_2(s\tau)\right]\right)\\
=&2s\left(C_1(s\tau)+sC_2(s\tau)\right).
\end{align*} 

\begin{remark}Bessel--Clifford functions are not so widely used since they can be written in terms the well-known modified Bessel functions $I_\nu$ by means of the equation  $\left(\frac{z}{2}\right)^{-\nu}I_\nu(z)=C_\nu(\frac{z^2}{4})$, see  \cite[p. 375]{AS}. This formula can help us to determine the growth of the continuous binomial coefficients: from the asymptotic expansion of $I_\nu(z)$ for fixed $\nu$ and large values of $z$ \cite[p. 377]{AS}, we find that $$C_\nu(z)\sim z^{-\left(\frac{\nu}{2}+\frac{1}{4}\right)} \frac{e^{2\sqrt{z}}}{2\sqrt{\pi}},\quad \text{ as }z\longrightarrow \infty \text{ and }
	|\text{arg}(z)|<\pi,$$ where the powers of $z$ assume their principal value (here $f(z)\sim g(z)$ means that $\lim_{z\rightarrow \infty} f(z)/g(z)=1$ in the corresponding domain). Then, it follows from equation (\ref{BC in terms BC}) that 
	$${\tau+s\brace s}\sim \left(\sqrt{s}+\sqrt{\tau}\right)^2 \frac{e^{2\sqrt{s\tau}}}{2\sqrt{\pi}(s\tau)^{3/4}},\quad \text{ as }s\tau\longrightarrow \infty,\quad |\text{arg}(s\tau)|<\pi.$$
\end{remark}

\

\section{Adding all the lengths}\label{Sec: Adding all the lengths}

When working with the Riemannian metric (\ref{Metric General}), we find that the paths $\phi_1(s)=\Psi(x_0+s,y_0)$ and $\phi_2(s)=\Psi(x_0,y_0+s)$, $0\leq s \leq s_0$, have lengths \begin{equation}\label{formulas length}l(\phi_1)=h(x_0+s_0)-h(x_0), \quad l(\phi_2)=f'(x_0)s_0,\end{equation} respectively. Regarding the Gaussian curvature, it only depends on the $x$--coordinate and it is given by the formula $-\frac{\left(f''/h'\right)'(x)}{h'(x)f'(x)}$, which can be easily checked by direct computations.

Our main example consists of surfaces of revolution. To fix ideas, consider an injective smooth curve  $\eta:(a,b)\rightarrow \R^3$, $\eta(s)=(\a(s),0,\beta(s))$ with $\a(s)>0$. When we rotate it around the $z$--axis, the map $\Psi(s,\phi)=(\a(s)\cos\phi,\alpha(s)\sin\phi,\beta(s))$, $(s,\phi)\in (a,b)\times (0,2\pi)$ provides a parameterization and the associated metric is given by  $$\left\|\eta'(s)\right\|^2ds^2+\a(s)^2d\phi^2, \quad  h(s)=\int_{s_0}^s \left\|\eta'(\tau)\right\|d\tau,\quad  f(s)=\int_{s_0}^s \alpha(\tau)d\tau.$$ Note that if $\eta$ is parameterized by arc-length, then $h(s)=s-s_0$. These surfaces 
realize the Riemannian metric (\ref{Metric General}) when $|f''(s)|< h'(s)$ because in this case we can simply take $\alpha(s)=f'(s)$ and $\beta'(s)=\pm\sqrt{h'(s)^2-f''(s)^2}$. 

If we choose $h(s)=s$ and $f(s)=-\cos(s)$, $0<s<\pi$, the previous conditions are fulfilled and we obtain a typical local parameterization of the sphere $S^2$. On the other hand, by taking $h(y)=f(y)=\ln(y)$, the metric would be $y^{-2}(dx^2+dy^2)$ which is the metric of the Poincar\'{e} plane $\mathbb{H}$. Thus, if we put $$\a(s)=\frac{1}{s},\quad \beta(s)=\pm\int_1^s \frac{1}{\tau}\sqrt{1-\frac{1}{\tau^2}}d\tau=\mp\sqrt{1-\frac{1}{s^2}}\pm\ln\left(s+\sqrt{s^2-1}\right),\quad s>1,$$ and we change $1/s=\sin(t)$, we recover the usual parameterization of the pseudosphere obtained by rotating the tractrix   $$\eta(t)=\left(\sin(t),0,\cos(t)+\ln\tan\left(\frac{t}{2}\right)\right).$$ This is a model of an open surface of constant Gaussian curvature equal to $-1$, see e.g., \cite[p. 198]{AbateTovena}.

Naturally, if $h(s)=f(s)=1$ we obtain the cylinder, which has curvature identically zero. Yet another model for curvature zero is the punctured plane in polar coordinates, i.e., $\eta(s)=(s,0,0)$, $s>0$, $h(s)=s, f(s)=s^2/2$, and the metric $ds^2+s^2d\phi^2$. In these cases Theorem \ref{Th. Main} shows that the integral of the length is given by the measure of $\Gamma_{p,q}(t)$ times the average of the lengths of the paths with configurations $(1,2)$ and $(2,1)$, since $f$ is quadratic.

\

We are now in position to define and compute the integral of the length over $\Gamma_{p,q}(t)$, and thus to prove Theorem \ref{Th. Main}. The first observation is that the length of the horizontal part of any path in $\Gamma_{p,q}(t)$ is actually constant, and it is given by $h(x_1)-h(x_0)$. Then, it is enough to calculate the integral of the length of the vertical parts. For this, we use the decomposition (\ref{Gamma pq}) and distinguish the four cases I.1, I.2, II.1 and II.2 as explained in Section \ref{Sec: Paths as simplexes and their measure}. To simplify notations we will write $a=x_1-x_0$ and $t-a=y_1-y_0$.

Let us start with the case I.1, where $\Gamma_{p,q}^{\textbf{c}}(t)$ is identified with $\Delta_{m-1}^a\times \Delta_{m-1}^{t-a}$, and thus we use the  variables $l_j=\sum_{k=0}^{j-1} s_{2k}$,  $\overline{l}_{j}=\sum_{k=0}^{j-1} s_{2k+1}$, $j=1,\dots,m-1,$ ($l_0=0$). By using the second formula in (\ref{formulas length}), we see that the length of the vertical parts of a path corresponding to $\textbf{s}\in \Delta_{m-1}^a\times \Delta_{m-1}^{t-a}$ is given by $$\sum_{j=0}^{m-1} f'\left(x_0+\sum_{i=0}^j s_{2i}\right)s_{2j+1}=f'(x_1)(t-a-\overline{l}_{m-1})+\sum_{j=0}^{m-2} f'\left(x_0+l_{j+1}\right)(\overline{l}_{j+1}-\overline{l}_j).$$ When we integrate this function over $W_{m-1}^a\times W_{m-1}^{t-a}$, we find that it is given by \begin{equation}\label{1 integral}
f'(x_1)\frac{a^{m-1}}{(m-1)!}\frac{(t-a)^m}{m!}+I_{m-1}(a,t-a),
\end{equation}
where $I_{m-1}(a,t-a)$ is equal to  \begin{align*}
&\int\limits_0^a\int\limits_0^{l_{m-1}}\cdots \int\limits_0^{l_2}\int\limits_0^{t-a}\int\limits_0^{\overline{l}_{m-1}}\cdots \int\limits_0^{\overline{l}_2} \sum_{j=0}^{m-2} f'\left(x_0+l_{j+1}\right)(\overline{l}_{j+1}-\overline{l}_j)d\overline{l}_1\cdots d\overline{l}_{m-1}dl_1\cdots dl_{m-1}\\
=&\int\limits_0^a\int\limits_0^{t-a} I_{m-2}(l_{m-1},\overline{l}_{m-1})d\overline{l}_{m-1}dl_{m-1}+\frac{(t-a)^{m}}{m!}\int\limits_0^{a}f'(x_0+l_{m-1})\frac{l_{m-1}^{m-2}}{(m-2)!}dl_{m-1}.
\end{align*}

To solve this recurrence and the ones that will appear in the other cases, we can use the following

\begin{lema}\label{Lema recurrencias} Let $f:\R\rightarrow\R$ be a differentiable function, and consider $r, r_1, r_2\in \N$, $\lambda, x_0\in\R$, and $a,b\in\R_{\geq0}$. If we define $I_0(a,b):=\lambda\frac{a^{r_1}}{r_1!}\frac{b^{r_2}}{r_2!}$ and $$ I_m(a,b):=\int\limits_0^a\int\limits_0^b I_{m-1}(x,y)dydx+\frac{b^{m+r}}{(m+r)!}\int\limits_0^a \frac{x^{m-1}}{(m-1)!}f'(x_0+x)dx,$$ then  $$I_m(a,b)=\lambda\frac{a^{r_1+m}}{(r_1+m)!}\frac{b^{r_2+m}}{(r_2+m)!}+\frac{a^{m-1}}{(m-1)!}\frac{b^{m+r}}{(m+r)!}\left(f(x_0+a)-f(x_0)\right),\quad m\in\N^*.$$
\end{lema}

\begin{proof}For $m=1$ the formula is clear, so we assume it is valid for some $m\in\N^\ast$. Using the definition and the induction hypothesis we see that \begin{align*}
	I_{m+1}(a&,b)=\int\limits_0^a\int\limits_0^b \lambda\frac{x^{r_1+m}}{(r_1+m)!}\frac{y^{r_2+m}}{(r_2+m)!}+\frac{x^{m-1}}{(m-1)!}\frac{y^{m+r}}{(m+r)!}\left(f(x_0+x)-f(x_0)\right)dydx\\
	&\qquad +\frac{b^{m+1+r}}{(m+1+r)!}\int\limits_0^a \frac{x^{m}}{m!}f'(x_0+x)dx\\
	=&\lambda\frac{a^{r_1+m+1}}{(r_1+m+1)!}\frac{b^{r_2+m+1}}{(r_2+m+1)!}+\frac{b^{m+1+r}}{(m+r+1)!}\int\limits_0^a \frac{d}{dx}\left(\frac{x^{m}}{m!}(f(x_0+x)-f(x_0))\right)dx\\
	=&\lambda\frac{a^{r_1+m+1}}{(r_1+m+1)!}\frac{b^{r_2+m+1}}{(r_2+m+1)!}+\frac{a^{m}}{m!}\frac{b^{m+1+r}}{(m+1+r)!}(f(x_0+a)-f(x_0)).
	\end{align*} The result follows from the principle of induction.
\end{proof}

By taking $\lambda=0$, $I_0=0$, and $r=1$ in Lemma \ref{Lema recurrencias},  we find that (\ref{1 integral}) becomes 
\begin{equation}\label{Integral 1}f'(x_1)\frac{a^{m-1}}{(m-1)!}\frac{(t-a)^{m}}{m!}+\left(f(x_1)-f(x_0)\right)\frac{a^{m-2}}{(m-2)!}\frac{(t-a)^{m}}{m!}.\end{equation} Note that any term with negative exponent is interpreted to be zero. 

We proceed now with case I.2, where $\Gamma_{p,q}^{\textbf{c}}(t)$ is again $\Delta_{m-1}^a\times \Delta_{m-1}^{t-a}$. In the variables $\overline{l}_j=\sum_{k=0}^{j-1} s_{2k}$, $l_j=\sum_{k=0}^{j-1} s_{2k+1}$, $j=1,\dots,m-1$, ($\overline{l}_m=t-a$) the length of the vertical parts of a path is given by
$$f'(x_0)\overline{l}_1+\sum_{j=1}^{m-1} f'(x_0+l_j)(\overline{l}_{j+1}-\overline{l}_j).$$ The integral of this function over $W_{m-1}^a\times W_{m-1}^{t-a}$ can be found as before by solving the same recurrence. Therefore, it is given by 
\begin{equation}\label{Integral 2}f'(x_0)\frac{a^{m-1}}{(m-1)!}\frac{(t-a)^m}{m!}+\left(f(x_1)-f(x_0)\right)\frac{a^{m-2}}{(m-2)!}\frac{(t-a)^m}{m!}.\end{equation}

For the case II.1,  $\Gamma_{p,q}^{\textbf{c}}(t)$ corresponds to $\Delta_{m}^a\times \Delta_{m-1}^{t-a}$ and we work with the variables $l_i=\sum_{k=0}^{i-1} s_{2k}$, $i=1,\dots,m$, $\overline{l}_j=\sum_{k=0}^{j-1} s_{2k+1}$, $j=1,\dots,m-1$ ($\overline{l}_0=0$ and  $\overline{l}_m=t-a$). The length of the vertical parts of a path in this case is $$\sum_{j=0}^{m-1} f'(x_0+l_{j+1})(\overline{l}_{j+1}-\overline{l}_j),$$ and the recurrence needed to calculate the integral of this function over  $W_{m}^a\times W_{m-1}^{t-a}$ is the same as in case (I.1) but with $r=0$. Thus, we get \begin{equation}\label{Integral 3}\left(f(x_1)-f(x_0)\right)\frac{a^{m-1}}{(m-1)!}\frac{(t-a)^m}{m!}.\end{equation}

Finally, for the case II.2, $\Gamma_{p,q}^{\textbf{c}}(t)$ corresponds to  $\Delta_{m-1}^a\times \Delta_{m}^{t-a}$, the variables are  $\overline{l}_i=\sum_{k=0}^{i-1} s_{2k}$, $i=1,\dots,m$, $l_j=\sum_{k=0}^{j-1} s_{2k+1}$, $j=1,\dots,m-1$, and the length of the vertical parts is given by
$$f'(x_0)\overline{l}_1+f'(x_1)(t-a-\overline{l}_{m})+\sum_{j=1}^{m-1} (\overline{l}_{j+1}-\overline{l}_j)f'(x_0+l_{j}).$$ The recurrence here for the integral of the last sum over  $W_{m-1}^a\times W_{m}^{t-a}$ is the same as in case (I.1) but with $r=2$. Thus, the integral is given by  \begin{equation}\label{Integral 4}
\frac{a^{m-1}}{(m-1)!}\frac{(t-a)^{m+1}}{(m+1)!}\left(f'(x_1)+f'(x_0)\right)
+\frac{a^{m-2}}{(m-2)!}\frac{(t-a)^{m+1}}{(m+1)!}\left(f(x_1)-f(x_0)\right).	
\end{equation}

Now, we define $\int_{\Gamma_{p,q}(t)} l(\gamma)d(\gamma)$, the \textit{integral of the length over} $\Gamma_{p,q}(t)$, as the sum of the integrals of the length of each of its part. For the horizontal parts, the result is $\left( h(x_1)-h(x_0) \right){t\brace x_1-x_0}$, being this length the constant $h(x_1)-h(x_0)$. For the vertical parts, we add together (\ref{Integral 1}), (\ref{Integral 2}), (\ref{Integral 3}) and (\ref{Integral 4}) to get 

$\left(f'(x_0)+f'(x_1)+f(x_1)-f(x_0)\right)\sum_{m=1}^\infty \frac{a^{m-1}}{(m-1)!}\frac{(t-a)^{m}}{m!}+2\left(f(x_1)-f(x_0)\right)\sum_{m=2}^\infty\frac{a^{m-2}}{(m-2)!}\frac{(t-a)^{m}}{m!}\linebreak+\left(f'(x_0)+f'(x_1)\right)\sum_{m=1}^\infty \frac{a^{m-1}}{(m-1)!}\frac{(t-a)^{m+1}}{(m+1)!}+\left(f(x_1)-f(x_0)\right)\sum_{m=2}^\infty \frac{a^{m-2}}{(m-2)!}\frac{(t-a)^{m+1}}{(m+1)!}$.

Then, by grouping common factors and taking into account formulas (\ref{Formula 1}) and (\ref{Formula 2}), the previous term is equal to $$\left(f(x_1)-f(x_0)\right) \frac{\d^2 F}{\d \tau^2}(a,t-a)+\frac{f'(x_1)+f'(x_0)}{2}\left((t-a)F(a,t-a)-a\frac{\d^2 F}{\d \tau^2}(a,t-a)\right).$$ This proves formula (\ref{Main formula}) in Theorem \ref{Th. Main}. Finally, for the last statement of the theorem, note that it is equivalent to show that $f$ satisfies the differential equation $$\frac{x_1-x_0}{2}(f'(x_1)+f'(x_0))=f(x_1)-f(x_0),$$ whose solutions are precisely quadratic functions.


\begin{thebibliography}{99} \small
	
	
	
	\bibitem{AbateTovena}{Abate M. and Tovena F.}, \emph{Curves and Surfaces}. Springer-Verlag, Milan, 2012.
	
	\bibitem{AS}{Abramowitz M. and Stegun I.},  \emph{Handbook of Mathematical Functions with Formulas, Graphs, and Mathematical Tables.} Dover, New York, ninth Dover printing, tenth GPO printing edition, 1964.
	
	\bibitem{CD1}{Cano L. and D\'{i}az R.}:
	``Indirect Influences on Directed Manifolds'',  \textit{Adv. Stud. Contemp. Math.}, 28 (2018) n° 1, 93--114.
	
	\bibitem{CD}{Cano L. and D\'{i}az R.}:
	Continuous Analogues for the Binomial Coefficients and the Catalan Numbers. Submitted to publication. Available at https://arxiv.org/pdf/1602.09132v4.pdf
	
	\bibitem{CD2001}{Chaichian M. and Demichev A.}:
	\emph{Path Integrals in Physics: Volume I Stochastic Processes and Quantum Mechanics}. Series in mathematical and computational physics. Taylor \& Francis, 2001.
	
	\bibitem{Feynman}{Feynman R. and Hibbs A.}: \emph{Quantum mechanics and integrals}. McGraw Hill, 1985.
	
	\bibitem{Salwinski}{Salwinski D.}: ``The Continuous Binomial Coefficient: An Elementary Approach'', \textit{Amer. Math. Monthly}, 123 (2018), no. 3, 231--244.
	
	\bibitem{W2}{Wakhare T. and Vignat C.}: ``A continuous analogue of lattice path enumeration: part II''. \textit{Online J. Anal. Comb}, Issue 14 (2019), $\# 04$.
	
\end{thebibliography}
\end{document}